\newif\ifpictures
\picturestrue

\documentclass[12pt]{amsart}

\usepackage{amssymb,amsmath}
\usepackage{hyperref}
\usepackage[dvips]{graphicx}

\ifpictures
\usepackage{psfrag}
\fi

\numberwithin{equation}{section}
\newtheorem{theorem}{Theorem}[section]
\newtheorem{lemma}[theorem]{Lemma}
\newtheorem{cor}[theorem]{Corollary}
\newtheorem{prop}[theorem]{Proposition}

\theoremstyle{definition}
\newtheorem{defini}[theorem]{Definition}
\newtheorem{example}[theorem]{Example}
\newtheorem{remark}[theorem]{Remark}
\newtheorem*{remark*}{Remark}

\title{Combinatorics and Genus of Tropical Intersections and Ehrhart Theory}
\author{Reinhard Steffens and Thorsten Theobald}

\thanks{Research supported by the Deutsche Forschungsgemeinschaft (DFG)}

\address{FB 12 -- Institut f\"ur Mathematik, Goethe-Universit\"at,
  Robert-Mayer-Str.\ 6--10, D--60325 Frankfurt am Main, Germany}

\email{steffens@math.uni-frankfurt.de, theobald@math.uni-frankfurt.de}
\newcommand{\C}{{\mathbb C}}
\newcommand{\R}{{\mathbb R}}

\newcommand{\Z}{{\mathbb Z}}
\newcommand{\N}{{\mathbb N}}
\newcommand{\Sph}{{\mathbb S}}

\newcommand{\vol}{\textup{vol}}
\newcommand{\mv}{\textup{MV}}

\DeclareMathOperator{\ME}{ME}
\DeclareMathOperator{\me}{me}

\begin{document}
\maketitle

\begin{abstract}
Let $g_1, \dots, g_k$ be tropical polynomials in $n$ variables 
with Newton polytopes $P_1, \dots, P_k$. We study combinatorial questions on the intersection of the tropical  hypersurfaces defined by $g_1, \dots, g_k$, such as the $f$-vector, the number of unbounded faces and (in case of a curve) the genus. 
Our point of departure is Vigeland's work \cite{Vigeland} who considered the special case
$k=n-1$ and where all Newton polytopes are standard simplices. 
We generalize these results to arbitrary $k$ and  arbitrary Newton polytopes $P_1, \dots, P_k$. 
This provides new formulas for the number of faces and the genus in terms of mixed volumes. By establishing
some aspects of a mixed version of Ehrhart theory we show that the genus of a tropical intersection curve equals the genus of a toric intersection curve corresponding to the same Newton polytopes.
\end{abstract}

\section{Introduction }
Tropical geometry allows to express certain algebraic-geometric problems in terms of discrete geometric
problems using \emph{correspondence theorems}. One general aim is to establish new
tropical methods to study the original algebraic problem (see, e.g. \cite{SturmfelsFeichtner,Draisma08, EinsiedlerKapranov, Mikhalkin04}). A prominent example is the work of Mikhalkin \cite{Mikhalkin04b} who gave a tropical formula for the number of plane curves of given degree and genus passing through a given number of points; see also \cite{GathmannMarkwig,Itenberg, SiebertNishinou} for related theorems. Providing methods to establish these correspondence statements is an important task of current research. 
Another important objective in tropical geometry is to understand the combinatorial structure
of the tropical varieties which can be regarded as polyhedral complexes in $n$-dimensional space.
See, e.g., \cite{SpeyerTropicalLinearSpaces, SpeyerSturmfels}.

In this paper, 
we consider a transversal intersection of tropical hypersurfaces given by polynomials $g_1, \ldots, g_k$ in $\R^n$ with Newton polytopes $P_1, \ldots,P_k$. For the special case $k = n-1$ and all $P_i$ standard simplices, Vigeland studied the number of vertices and unbounded edges as well as the genus of this intersection \cite{Vigeland}.
His methods strongly rely on the special structure of the Newton polytopes.

The new contributions of the present paper can be stated as follows.
Firstly, we provide a uniform and systematic treatment of the whole $f$-vector (i.e., the vector of face numbers) of the tropical transversal intersection. In particular, we show how to reduce these counts to well-established tropical intersection theorems. Generalizing the results in \cite{Vigeland}, our approach also covers the general mixed case, where we start from polynomials $g_1, \ldots, g_k$ with arbitrary Newton polytopes $P_1, \ldots, P_k$. We obtain formulas expressing the number of faces (Theorems \ref{theorem:VertexCount} and \ref{theorem:unbounded}) and the genus (Theorem \ref{Theorem:Genus}) in terms of mixed volumes.

Secondly, we establish a combinatorial connection from the tropical genus of a curve to the genus of a toric curve corresponding to the same Newton polytopes. In \cite{Khovanskii78}, Khovanskii gave a characterization of the genus of a toric variety in terms of integer points in Minkowski sums of polytopes. 
We show that in the case of a curve this toric genus coincides with the tropical genus (Theorem \ref{the:CompareGenus}). More than this result itself, we think that the methods to establish it
are of particular interest. Khovanskii's formula applies integer points in Minkowski sums of polytopes, whereas the mentioned formula for the tropical genus is given in terms
of mixed volumes. For the special case $n=2$ the connection boils down to the classical Theorem of Pick relating the number of integer points in a polygon to its area. We develop a Pick-type formula for the surface volume of a lattice complex in terms of integer points (Theorem \ref{the:PickOnSurface}) to show that in the generalized \emph{unmixed} case ($n$ arbitrary, all $P_i$ coincide) the connection reduces to certain $n$-dimensional generalizations of Pick's theorem (Macdonald \cite{MacDonald62}). To approach the general mixed case we have to develop new aspects of a mixed Ehrhart theory  (Theorem \ref{the:MixedEhrhart}) based on McMullen's results on valuations \cite{McMullen77}.

The paper is structured as follows. Section \ref{se:preliminaries} gives a short introduction to tropical geometry, mixed volumes and Ehrhart theory. In Section \ref{section:f-vector} we study the number of bounded and unbounded faces in a tropical intersection. Using these results we give in Section \ref{section:LatticePointsAndGenus} a formula for the genus of a tropical intersection curve and 
show that this tropical genus coincides with the toric genus that depends on the same Newton polytopes.

\section{Preliminaries\label{se:preliminaries}}
\subsection{Tropical geometry} 
Let $\R_{\textrm{trop}}:=(\R\cup\{-\infty\},\oplus,\odot)$ denote the
\emph{tropical semiring}. The arithmetic operations of \emph{tropical addition} $\oplus$ and \emph{tropical multiplication} $\odot$ are
\[
x \oplus y = \text{max}\{x,y\} \quad \text{and} \quad x \odot y = x+y \ .
\]
Equivalently tropical addition can be defined as $\text{min}\{x,y\}$ (e.g. \cite{FirstStepsTropical}) but results in either preferred notation can easily be translated into each other.
A \emph{tropical Laurent polynomial} $f$ in $n$ variables $x_1,\dots,x_n$ is an expression of the form
\begin{equation}  \label{Def:TropPoly}
f= \bigoplus_{\alpha \in \mathcal{S}(f)} c_{\alpha} \odot x_1^{\alpha_1} \odot \dots \odot x_n^{\alpha_n} = \textup{max}_{\alpha \in \mathcal{S}(f)} (c_{\alpha} + \alpha_1x_1 + \dots + \alpha_nx_n)  
\end{equation}
with real numbers $c_\alpha$. The support set $\mathcal{S}(f)$ is always assumed to be a finite subset of $\Z^n$, and its convex hull $P(f)\subset\R^n$ is called the \emph{Newton polytope of $f$}. A tropical polynomial $f(x_1,\dots,x_n)$ defines a convex, piecewise linear function $f:\R^n \rightarrow \R$ and we define the \emph{tropical hypersurface} $X(f)$ as the non-linear locus of $f$. See \cite{Mikhalkin04} or \cite{FirstStepsTropical} for more detailed introductions. 

Note that any tropical hypersurface $X(f)$ can be interpreted as a polyhedral complex of codimension $1$ in $\R^n$. The $m$-dimensional cells of a polyhedral complex $X$ will be denoted by $X^{(m)}$. For tropical polynomials $f_1, f_2$ we have $P(f_1 \odot f_2) = P(f_1)+P(f_2)$ and $X(f_1 \odot f_2) = X(f_1) \cup X(f_2)$.
\begin{example}\label{Example1}
Consider the two tropical polynomials
\begin{eqnarray*}
f&=& -62x \oplus 97x^{2} \oplus -73y^{2} \oplus -4x^{3}y \oplus -83x^{2}y^{2} \oplus -10y^{4} \ , \\
 g&=& -10x^{2}y \oplus 31x^{3}y \oplus -51xy^{3} \oplus 77y^{4} \oplus 95x^{2}y^{3} \oplus y^{5} \, ,
\end{eqnarray*}
where the tropical multiplication symbol $\odot$ is omitted.
\begin{figure}[bt]
\begin{center}
\includegraphics[width=2.8cm]{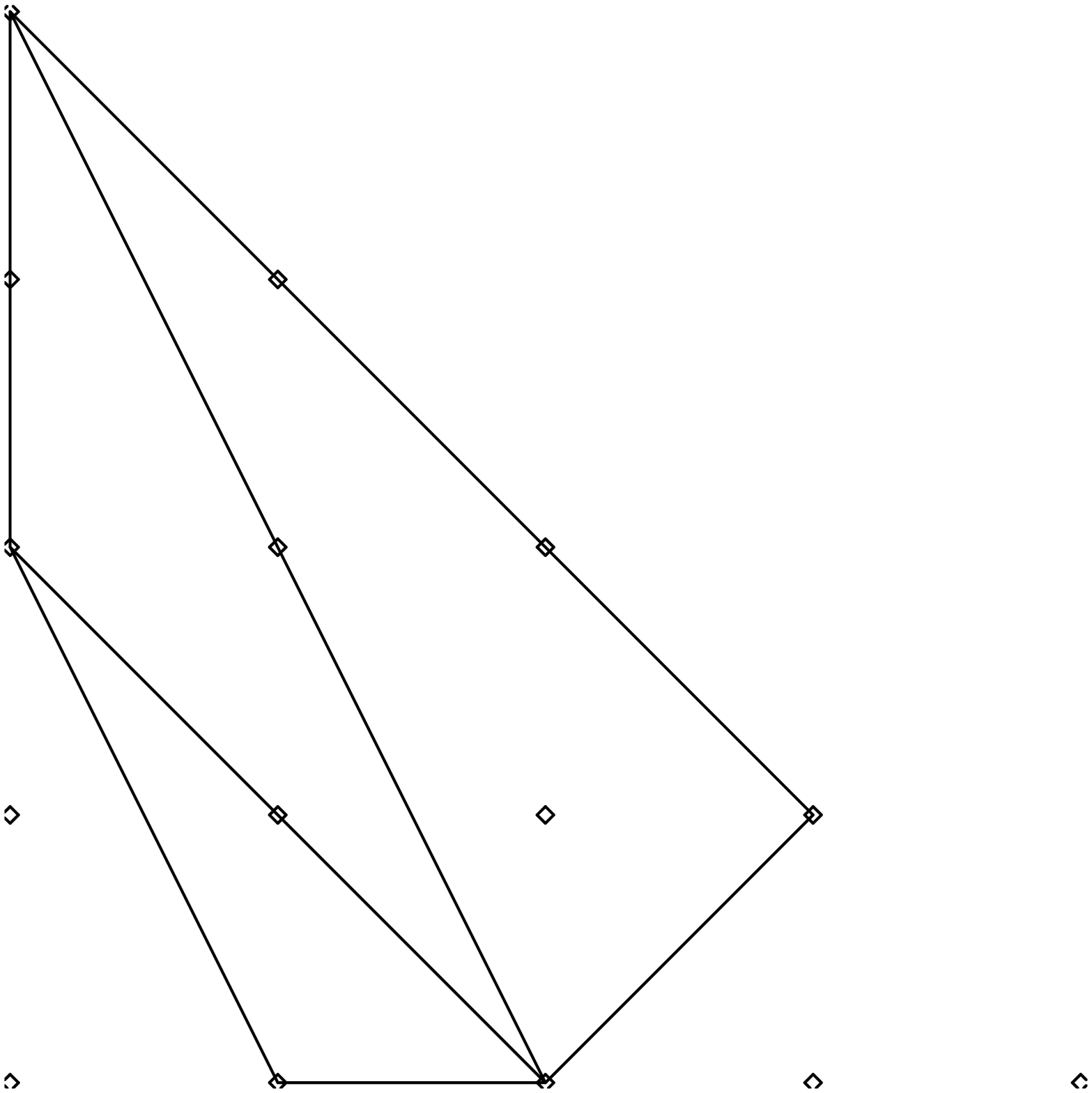}
\quad \includegraphics[width=3.2cm]{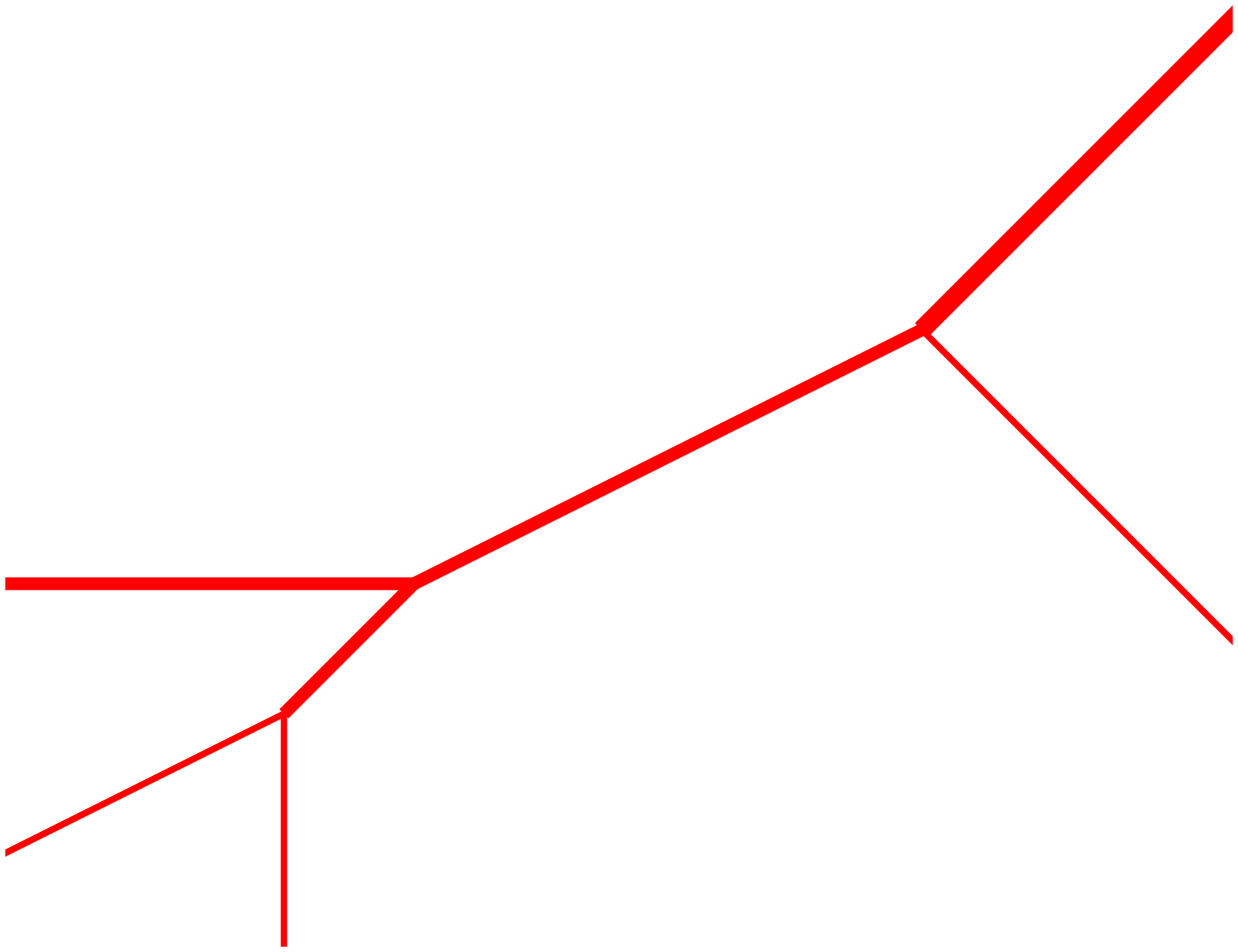}
\qquad \qquad \includegraphics[width=2.8cm]{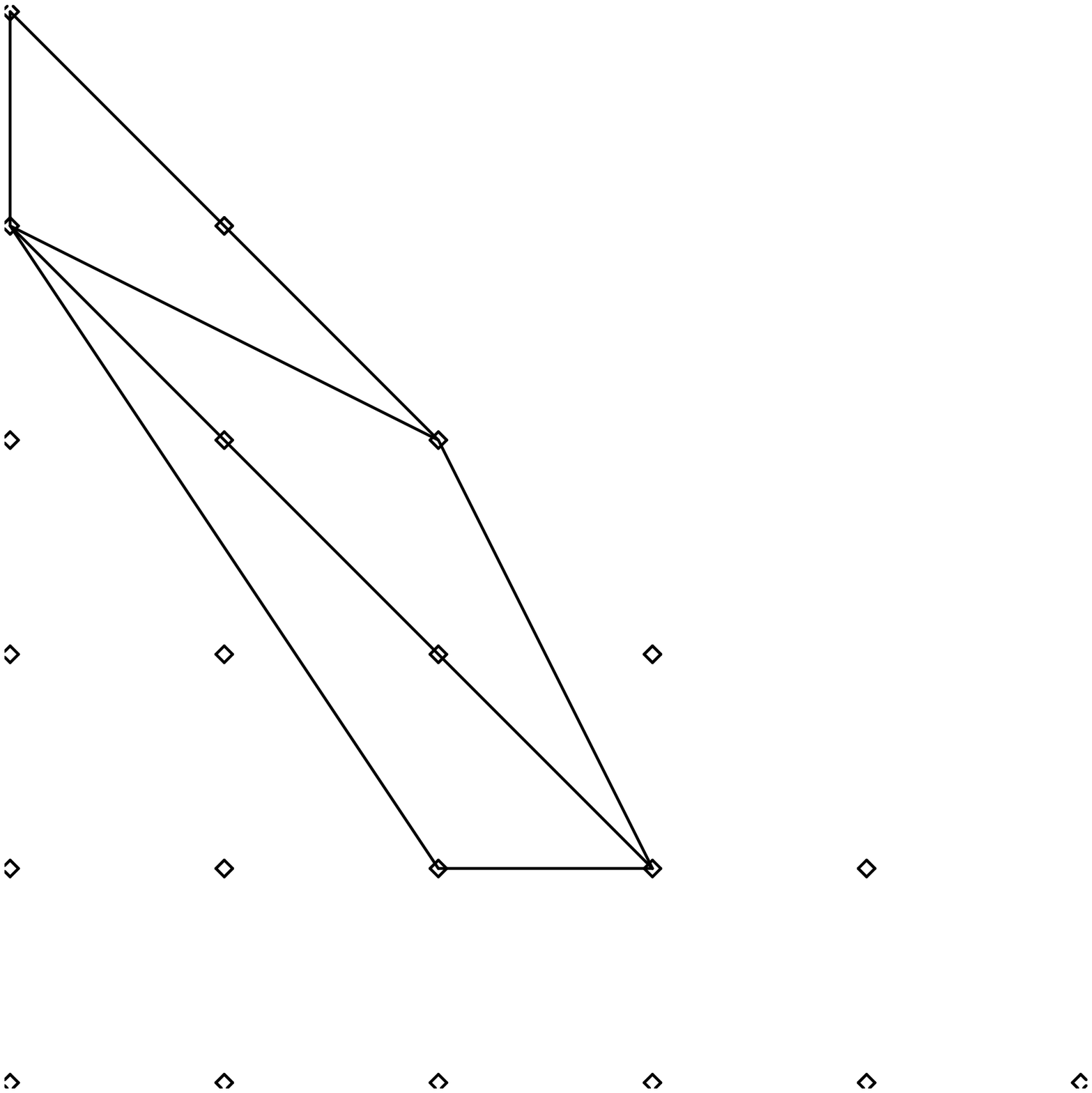}
\ \includegraphics[width=3.2cm]{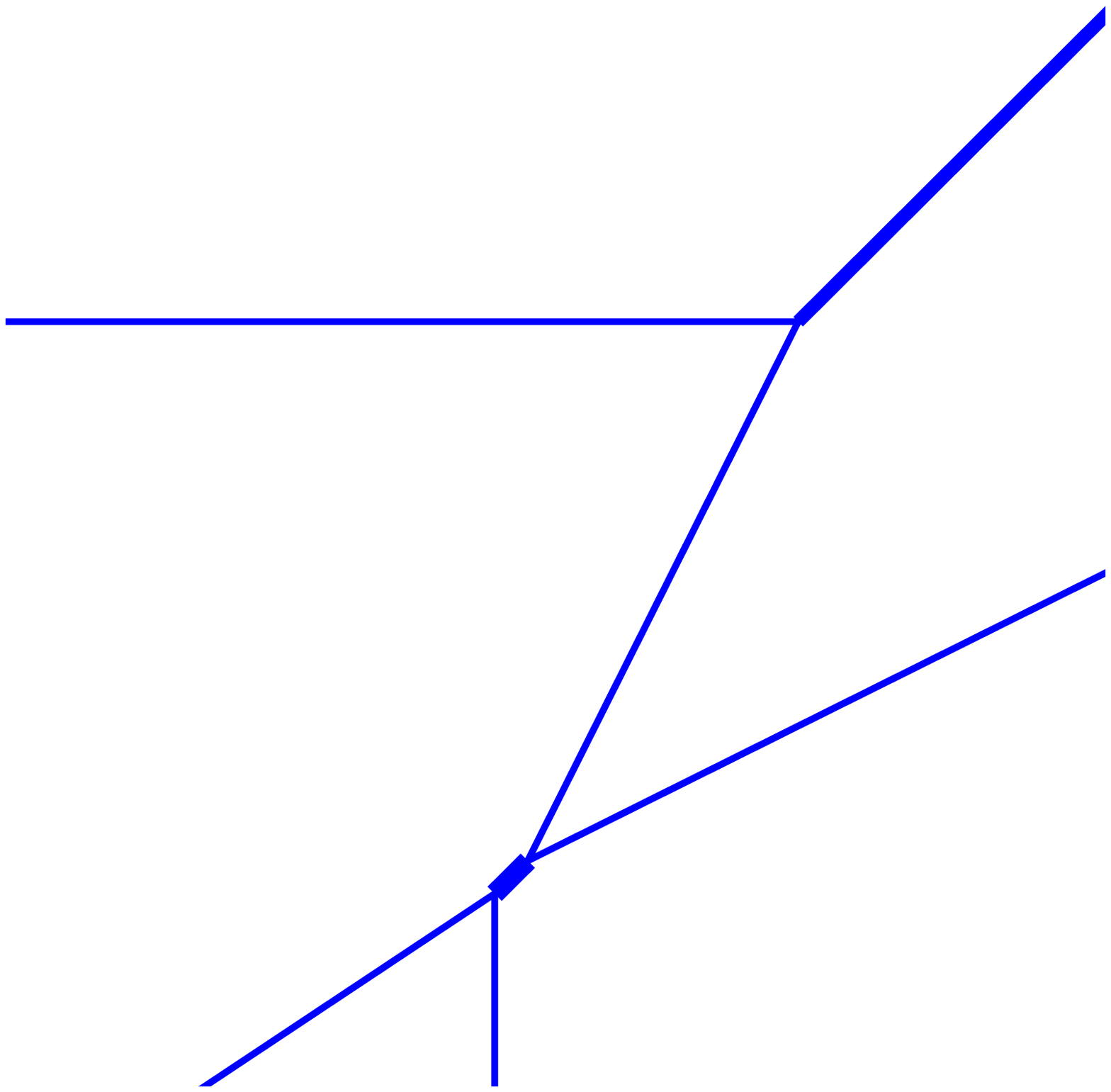}
\end{center}
\medskip
\begin{center}
\includegraphics[width=3.7cm]{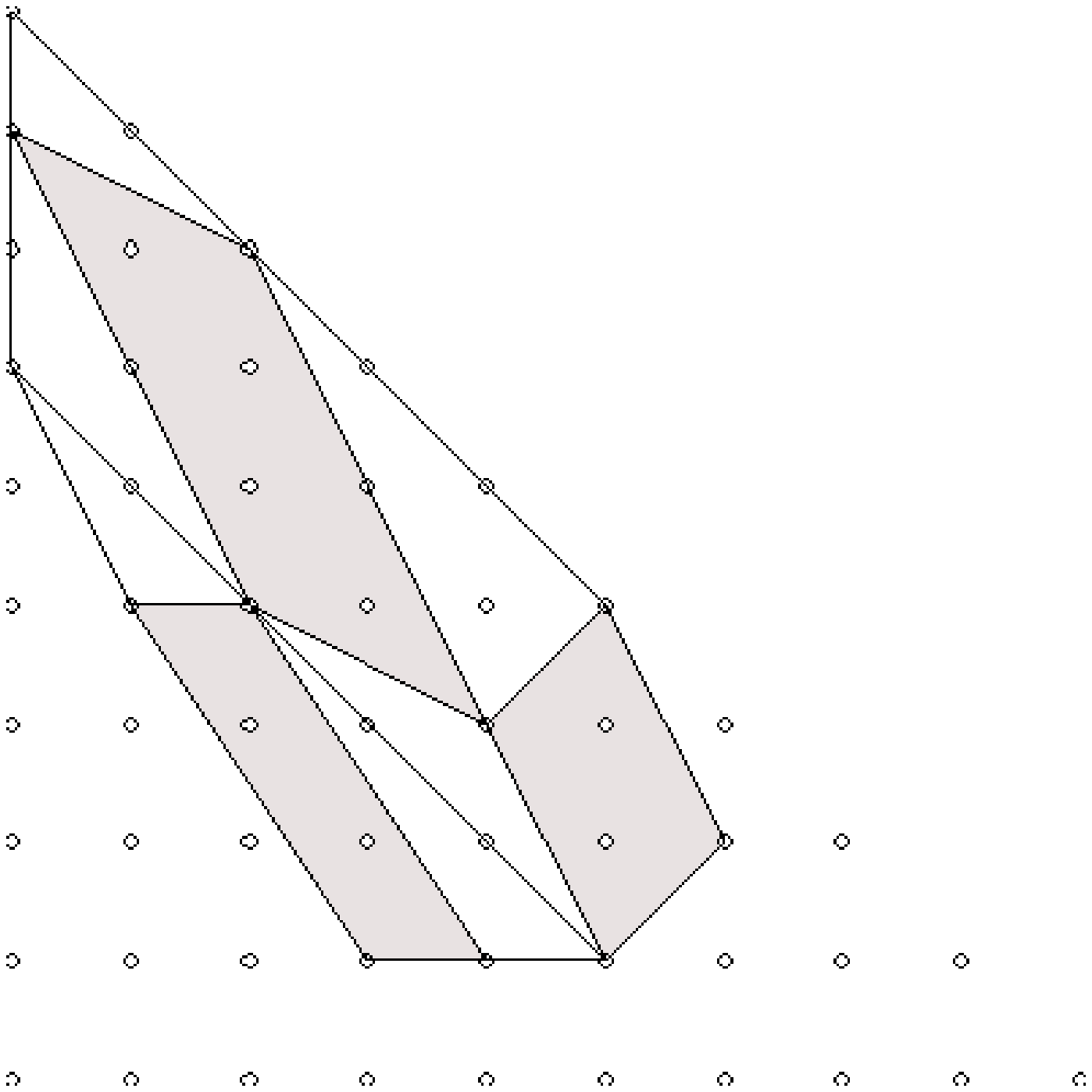}
\quad \includegraphics[width=3.7cm]{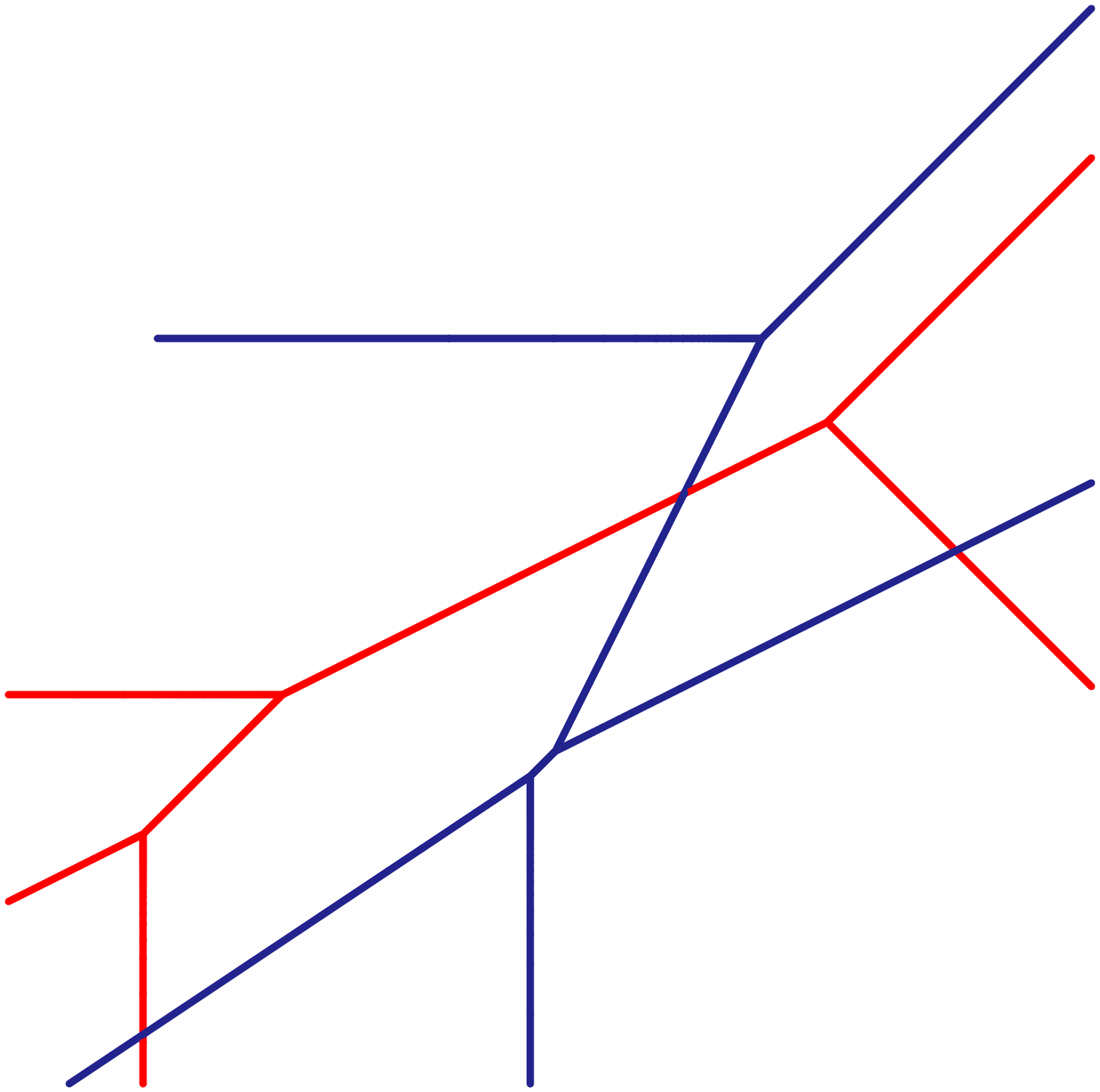}
\end{center}
\caption{Top: $P(f)$ and $P(g)$ with the privileged subdivision and the tropical curves $X(f)$ and $X(g)$. (Here bold edges indicate higher multiplicities.) Bottom: $P(f\odot g)$ with the privileged subdivision and the tropical curve $X(f\odot g)$.}
\label{Figure:f_and_g}
\end{figure}
Figure \ref{Figure:f_and_g} shows  their curves and their Newton polytopes as well as the Newton polytope of the product $f\odot g$ and the union $X(f)\cup X(g)$.
\end{example}

\subsection{Mixed volumes}
The \emph{Minkowski sum} of two point sets $A_1, A_2 \subset \R^n$ is defined as
\[
A_1+A_2 = \left\{a_1+a_2 \, |\, a_1\in A_1, a_2\in A_2 \right\} \ . 
\]
Let $P_1,\dots,P_n$ be $n$ polytopes in $\R^n$. For non-negative parameters 
$\lambda_1,\dots,\lambda_n$ the function
$\text{vol}_n(\lambda_1P_1+\dots+\lambda_nP_n)$ is a homogeneous
polynomial of degree $n$ in $\lambda_1,\dots,\lambda_n$ with
non-negative coefficients (see e.g. \cite{Schneider}).  The coefficient of the mixed
monomial $\lambda_1\cdots\lambda_n$ is called the {\it mixed volume of $P_1,\dots,P_n$} and is
denoted by $\text{MV}_n(P_1,\dots,P_n)$. 

We denote by $\textup{MV}_{n}(P_1,d_1; \dots; P_k,d_k)$ the mixed volume where $P_i$ is taken $d_i$ times and $\sum_{i=1}^k d_i =n$.
The mixed volume is linear in each argument, i.e.
\begin{equation}
\text{MV}_n(\dots, \alpha P_i + \beta P'_i,\dots)= 
\alpha\, \text{MV}_n(\dots, P_i ,\dots)+ \beta \, \text{MV}_n(\dots, P'_i ,\dots) \label{linearity} 
\end{equation}
and it generalizes the usual volume in the sense that
$
\text{MV}_n(P,\dots,P)=  n!\,\text{vol}_n (P)
$
holds (see \cite{Schneider}).

Let $P=P_1+\dots+P_k\subset\R^n$ be a Minkowski sum of polytopes that affinely spans $\R^n$. A sum $C=F_1+\dots+F_k$ of faces $F_i\subset P_i$ is called {\it cell} of $P$. A {\it subdivision} of $P$ is a collection $\Gamma=\{C_1,\dots,C_m\}$ of cells such that each cell is of full dimension, the intersection of two cells is a face of both and the union of all cells covers $P$. Each cell is given a type $type(C)= (\text{dim}(F_1),\dots, \text{dim}(F_k))$. Clearly the entries in the type vector sum up to at least the dimension of the cell $C$.  A subdivision is called {\it mixed} if for each cell $C\in \Gamma$ we have that $\sum d_i =n$ where $type(C)=(d_1,\dots,d_k)$. 
Cells of type $(d_1,\dots,d_k)$ with $d_i\geq 1$ for each $i$ will be called {\it mixed cells}.

With this terminology the mixed volume can be calculated by
\begin{equation}\label{MV_explicit}
\text{MV}_n(P_1,d_1;\dots;P_k,d_k) =\sum_{C} d_1!\, \cdots d_k!\ \text{vol}_n\,(C) 
\end{equation}
where the sum is over all cells $C$ of type $(d_1,\dots,d_k)$ in an arbitrary mixed subdivision of $P_1+\dots+P_k$ (see \cite{Huber}).

\subsection{Privileged subdivisions and duality}\label{section:duality}

The Newton polytope $P(f)$ of a tropical polynomial $f$ comes with a \emph{privileged subdivision} $\Gamma(f)$.
Namely we lift the points $\alpha\in\mathcal{S}(f) $ into $\R^{n+1}$
using the coefficients $c_{\alpha}$ as lifting values. 
The set of those facets of $\hat{P}(f):= \text{conv}\{(\alpha,c_{\alpha})\, | \, \alpha \in \mathcal{S}(f)\}$ which have an inward pointing normal with a negative last coordinate is called the \emph{upper hull}. If we project down this upper hull back to $\R^n$ by forgetting the last coordinate we get a subdivision of $P(f)$.  
On a set of $k$ tropical polynomials $f_1,\dots,f_k$ the coefficients induce a \emph{privileged subdivision} $\Gamma(f_1,\dots,f_k)$  of $P(f_1)+\dots+P(f_k)$ by projecting down the upper hull of $\hat{P}(f_1)+\dots+\hat{P}(f_k)$. For a generic choice of coefficients in the system $f_1,\dots,f_k$ this subdivision will be mixed (see \cite{Huber}).

The subdivision $\Gamma(f_1,\dots,f_k)$ and the union $X(f_1)\cup \dots \cup X(f_k)$ of tropical hypersurfaces are polyhedral complexes which are dual in the sense that there is a one-to-one correspondence  between their cells which reverses the inclusion relations (see \cite{BertrandBihan,MikhalkinPairsofPants}). Each cell $C$ in $\Gamma(f_1,\dots,f_k)$ corresponds to a cell $A$ in $X(f_1)\cup \dots \cup X(f_k)$ such that $\dim(C)+\dim(A)=n$, $C$ and $A$ span orthogonal real affine spaces and $A$ is unbounded if and only if $C$ lies on the boundary of $P(f_1)+\dots+P(f_k)$. Furthermore we have that a cell $A$ of $X(f_1)\cup \dots \cup X(f_k)$ is in the intersection $\mathcal{I}=X(f_1)\cap \dots \cap X(f_k)$ if and only if the corresponding dual cell $C$ in $\Gamma(f_1,\dots,f_k)$ is mixed.

A cell $A$ in $\mathcal{I}$ can be written as $A=\bigcap_{i=1}^k A_i$ where $A_i\in X_i$. If we require that $A$ lies in the relative interior of each $A_i$ then this representation is unique. The dual cell $C$ of $A$ has then a unique decomposition into a Minkowski sum $C=F_1+\dots+F_k$ where each $F_i$ is dual to $A_i$. We will always refer to this decomposition if not stated otherwise.

\subsection{Ehrhart theory and valuations}\label{section:Ehrhart}
Next we review some results on Ehrhart polynomials (see \cite{Barvinokbook, Ehrhart67})  and on valuations (see \cite{McMullen77}) of \emph{lattice polytopes}, i.e. polytopes with vertices in $\Z^n$. Let $B(P)$ and $B^{+}(P)$ denote the number of integer points and the number of interior integer points of a lattice polytope $P$, respectively.

Ehrhart showed that the number of integer points in $t\cdot P \cap \Z^n$ for $t\in\N$ is a polynomial in $t$, i.e.
\[
B(t\cdot P) = E_P(t)\quad  \text{ for some polynomial } E_P(x)=\sum_{i=0}^n e_i(P)\cdot x^i\ .
\]
The polynomial $E_P(t)$ is called the \emph{Ehrhart polynomial} of $P$ and its coefficients $e_i(P)$  are called \emph{Ehrhart coefficients}. The following identities are known for the coefficients:
\begin{equation}\label{eq:CoeffIdentities}
e_n(P)=\vol_n(P), \quad e_{n-1}(P)= \frac{1}{2}\sum_{F \text{ facet of }P}\vol_{n-1}'(F), \quad e_0(P)=1\ .
\end{equation}
For the remaining coefficients we do not have explicit expressions but via 
Morelli's formula the coefficient $e_k(P)$ can be expressed
in terms of the $k$-dimensional faces of $P$ (see \cite{McMullen93}).

Concerning the number of interior integer points similar results can be stated using the \emph{reciprocity law}
\begin{equation}\label{eq:EhrhartReciprocity}
B^+(t\cdot P)= (-1)^n E_P(-t) = (-1)^n \sum_{i=0}^n (-1)^i e_i(P)\cdot t^i\ .
\end{equation}

Let $\mathcal{P}(\Z^n)$ denote the space of all lattice polytopes in $\R^n$. A mapping $\varphi : \mathcal{P}(\Z^n) \rightarrow \R$ is called a \emph{$\Z^n$-valuation} if 
\[
\varphi(P\cup Q)+\varphi(P\cap Q) = \varphi(P)+ \varphi(Q)
\]
for $P\cup Q, P\cap Q, P, Q \in \mathcal{P}(\Z^n)$ and $\varphi(P+a)= \varphi(P)$ for all $a\in\Z^n$. If furthermore we have that 
\[
\varphi(t\cdot P) = t^r \varphi(P)
\]
holds, then the $\Z^n$-valuation $\varphi$ will be called \emph{homogeneous of degree $r$}.
Examples of $\Z^n$-valuations include the volume and the number of integer points in a lattice polytope.
The volume is a homogeneous valuation of degree $n$.

 For results in later sections it is crucial that the Ehrhart coefficients $e_r(P)$ are 
 $\Z^n$-valuations which are homogeneous of degree $r$ (see \cite{Betke85, McMullen77}). A key ingredient will be the following lemma by McMullen (see \cite{McMullen77}).
\begin{prop}[McMullen \cite{McMullen77}] \label{lemma:McMullen}
Let $\varphi_r$ be a homogeneous $\Z^n$-valuation of degree $r$ and let $t_1,\dots,t_k$ be integers. Then for any polytopes $P_1,\dots,P_k \in \mathcal{P}(\Z^n)$ we have
\begin{equation*}
\varphi_r(t_1\cdot P_1+\dots+t_k\cdot P_k) = \sum_{r_1,..,r_k}\binom{r}{r_1 \dots r_k} \varphi_r'( P_1,r_1;\dots;P_k,r_k) t_1^{r_1}\cdots t_k^{r_k}\ ,
\end{equation*}
where the multinomial coefficient $\binom{r}{r_1 \dots r_k}$ is defined by
\[
\binom{r}{r_1 \dots r_k}= \begin{cases}\frac{r!}{r_1!\cdots r_k!} & \text{if }r_i\geq0 \text{ and } \sum_i r_i = r\\ 0 & \text{ otherwise.}    \end{cases}
\]
\end{prop}
The coefficients $\varphi_r'( P_1,r_1;\dots;P_k,r_k)$ are called \emph{mixed $\Z^n$-valuations}. One can show that
\[
\varphi_r(P) = \varphi_r'(P,r) = \varphi_r'(\underbrace{P,\dots,P}_{r-\text{times}})
\]
and that $\varphi_r'$ is independent of $P_i$ if $r_i=0$ but we will not need this or any other explicit expression of these coefficients.
\begin{remark}
The results in \cite{McMullen77} are stated more generally for any additive subgroup $\Lambda$ of the Euclidean space. However here we just need the case of $\Z^n$. 
\end{remark}

\section{The $f$-vector of a tropical intersection}\label{section:f-vector}
Let $g_1,\dots,g_k$ be tropical polynomials in $n$ variables $x_1,\dots,x_n$ with Newton polytopes $P_1,\dots,P_k$ and let $X_i:=X(g_i)$ denote their tropical hypersurfaces in $\R^n$. 
After introducing some natural intersection multiplicities, we  study the number $f_j$ of $j$-dimensional faces in the polyhedral complex $\mathcal{I}=X_1\cap\dots\cap X_k$. The vector $(f_0,\dots,f_{n})$ is called the \emph{$f$-vector} of $\mathcal{I}$. The number of unbounded faces will be treated separately afterwards.

\subsection{Intersections of tropical hypersurfaces}
An intersection $\mathcal{I}=X_1\cap \dots\cap X_k$ is called \emph{proper} if $\textup{dim}(\mathcal{I})=n-k$. 
 $\mathcal{I}$ is \emph{transversal along a cell} $A$ of this complex if the dual cell $C=F_1+\dots+F_k$ in the privileged subdivision of $P_1+\dots+P_k$ satisfies
\[
\textup{dim}(C)= \textup{dim}(F_1)+\dots+\textup{dim}(F_k)\ .
\]
We call the intersection \emph{transversal} if for each subset $J\subset\{1,\dots,k\}$ the intersection is proper and transversal along each cell of the complex.  
In the dual picture a transversal intersection implies that the privileged subdivision of $P_1+\dots+P_k$ is mixed. Note that in a transversal intersection each cell $A$ of $\mathcal{I}$ lies in the relative interior of each cell $A_i$ from $X_i$ that is involved in the intersection.  

In the case of a non-transversal intersection $\mathcal{I}$ we can perturb the hypersurfaces by a small parameter $\varepsilon$ to obtain again a transversal intersection $\mathcal{I}_{\varepsilon}$. The \emph{stable intersection} $\mathcal{I}_{\textup{st}}$ is defined as the limit of these transversal intersections when $\varepsilon$ goes to $0$, 
\[
\mathcal{I}_{\textup{st}} = X_1 \cap_{ \textup{st}}\dots \cap_{ \textup{st}} X_k = \lim_{\varepsilon \rightarrow 0} X_1^{(\varepsilon_1)} \cap \dots \cap X_k^{(\varepsilon_k)} 
 \]
(see \cite{FirstStepsTropical}). Stable intersections are always proper and they have some more comfortable features. As mentioned above a tropical hypersurface $X(g)\subset\R^n$ is a polyhedral complex of dimension $n-1$. The stable intersection of $X(g)$ with itself gives the $(n-2)$-skeleton of $X(g)$.
In particular we can isolate the vertices of $X(g)$ by stably intersecting $X(g)$ $(n-1)$-times with itself.

Every face of a tropical intersection $\mathcal{I}$ naturally comes with a multiplicity.
We follow the notation of Bertrand and Bihan \cite{BertrandBihan}, whose approach is consistent
with those in \cite{Katz08, Mikhalkin04, sturmfels-tevelev-2008}.

\begin{defini}[Intersection multiplicity] \label{Def:Multiplicity}
Each cell $A$ in an intersection $\mathcal{I}$ can be assigned a \emph{multiplicity} (or \emph{weight}) as follows. Let $C=F_1+\dots+F_k$ be its dual cell in $P_1+\dots+P_k$.
 If $A$ is of dimension $j$ then $C$ is of dimension $n-j$ and we denote its type by $(d_1,\dots,d_k)$. 
For a transversal intersection define
\begin{align}
m_A &:= \left(\prod_{i=1}^k d_i!\cdot\vol'_{d_i}(F_i)\right)\cdot \vol'_{n-j}(\mathcal{P})
\nonumber\\
&  = \textup{MV}'_{n-j}(F_1,d_1;\dots;F_k,d_k)\label{Eq:Multiplicity}
\end{align} 
where $\mathcal{P}$ is a fundamental lattice polytope in the $(n-j)$-dimensional sublattice $\Z(F_1)+\dots+\Z(F_k)$ and where $\vol'_{d_i}$ denotes the volume in the lattice $\Z(F_i)$ spanned by the integer vectors of $F_i$.
(For more background on these relative volume forms and the proof that equality holds in \eqref{Eq:Multiplicity} see \cite{BertrandBihan}.)

In the non-transversal case we have that $n-j \leq d_1+\dots +d_k$ and we define,
\[
m_A := \sum_{\begin{smallmatrix}(e_1,\dots,e_k) \textup{ s.t.}\\ \sum e_i = n-j;\ 
e_i\leq d_i  \end{smallmatrix}}
\textup{MV}'_{n-j}(F_1,e_1;\dots;F_k,e_k)\ .
\]
\end{defini}

\begin{theorem}[Tropical Bernstein, see \cite{BertrandBihan, FirstStepsTropical}]\label{theorem:TropicalBernstein}
Suppose the tropical hypersurfaces $X_1 ,$ $ \dots , X_n \subset \R^n $ with Newton polytopes  $ P_1 , \dots , P_n $ intersect in finitely many points. Then the number of intersection points counted with multiplicity is  $\textup{MV}_n(P_1,\dots,P_n)$.

Furthermore the stable intersection of $n$ tropical hypersurfaces $X_1 ,$ $ \dots , X_n$ always consists of $\textup{MV}_n(P_1,\dots,P_n)$ points counted with multiplicities.
\end{theorem}

\subsection{The number of $j$-faces in $\mathcal{I}$}
Let $\mathcal{I}= X_1\cap\dots\cap X_k$ be a transversal intersection. Hence the intersection is proper which implies that the number of $j$-dimensional faces in $\mathcal{I}$ is $0$ if $j\geq n-k$. By using the duality approach described in Section \ref{section:duality} the number of $j$-faces can be expressed in terms of mixed volumes.
\begin{theorem}\label{theorem:VertexCount}
The number of $j$-faces in $\mathcal{I}$ counting multiplicities is
\begin{equation} \label{VertexCount}
\sum_{A \in \mathcal{I}^{(j)}} m_A = \sum_{\begin{smallmatrix}(d_1,\dots,d_k)\textup{ s.t.}\\ d_i\geq 1\textup{ and } \sum_i d_i = n-j   \end{smallmatrix}}  \textup{MV}'_{n-j} (P_1,d_1;\dots;P_k,d_k)   \ , 
\end{equation}
where $\textup{MV}'_{n-j} (P_1,d_1;\dots;P_k,d_k)$ is interpreted as in \eqref{MV_explicit} as the sum over the lattice volume times $d_1!\cdots d_k!$ of all  $(n-j)$-dimensional cells of type $(d_1,\dots,d_k)$ in a mixed subdivision of $P_1+\dots+P_k$.
\end{theorem}
Note that this implies the tropical version of Bernstein's theorem (see Theorem \ref{theorem:TropicalBernstein}) for $k=n$ and $j=0$.
\begin{proof}
Each $j$-dimensional cell $C$ in the mixed subdivision of $P_1+\dots+P_k$ is dual to an $(n-j)$-dimensional cell $A$ of $X_1\cup\dots\cup X_k$. If $C$ is a mixed cell, i.e. $d_i\geq 1$ for all $i$, its dual $A$ is contained in every $X_i$. Hence, by Definition \ref{Def:Multiplicity}
\[
\sum_{A \in \mathcal{I}^{(j)}} m_A = \sum_{\begin{smallmatrix}(d_1,\dots,d_k)\textup{ s.t.}\\ d_i\geq 1\textup{ and } \sum_i d_i = n-j   \end{smallmatrix}} \sum_{C=F_1+\dots+F_k} \mv'_{n-j}(F_1,d_1;\dots;F_k,d_k)
\]
where the second sum runs over all cells $C$ of type $(d_1,\dots,d_k)$. 
If we denote by $\vol_{d_i}'(F_i)$ the volume of $F_i$ in the lattice spanned by the integer points of $F_i$ and furthermore denote by $\mathcal{P}$ the fundamental lattice parallelotope in $\Z^{n-j}$ defined by $F_1,\dots,F_k$ then \eqref{Eq:Multiplicity} implies
\begin{align*}
\mv'_{n-j}(F_1,d_1;\dots;F_k,d_k) &= d_1!\cdots d_k!\, \vol'_{d_1}(F_1)\cdots \vol'_{d_k}(F_k)\, \vol'_{n-j}(\mathcal{P})\\
&= d_1!\cdots d_k!\, \vol'_{n-j}(C) \ .
\end{align*} 
Hence we have
\begin{align*}
\sum_{A \in \mathcal{I}^{(j)}} m_A &= \sum_{\begin{smallmatrix}(d_1,\dots,d_k)\text{ s.t.}\\ d_i\geq 1\text{ and } \sum_i d_i = n-j   \end{smallmatrix}} \sum_{\begin{smallmatrix} C \text{ of type}\\ (d_1,\dots,d_k) \end{smallmatrix} } d_1!\cdots d_k!\, \vol'_{n-j}(C)\\
&= \sum_{\begin{smallmatrix}(d_1,\dots,d_k)\text{ s.t.}\\ d_i\geq 1\text{ and } \sum_i d_i = n-j   \end{smallmatrix}}  \textup{MV}'_{n-j} (P_1,d_1;\dots;P_k,d_k)
\end{align*} 
where we used \eqref{MV_explicit} for the last identity.
\end{proof}

In Section \ref{section:LatticePointsAndGenus} we focus on the number of  vertices in tropical intersection curves. Hence we state Theorem \ref{theorem:VertexCount} for $k=n-1$ and $j=0$ again which gives a much nicer expression. 
\begin{cor}\label{Cor:TransvIntVertexCount}
Let $\mathcal{I}=X_1 \cap \dots \cap X_{n-1}$ be a transversal intersection curve in $\R^n$ of $n-1$ tropical hypersurfaces with corresponding Newton polytopes $P_1,\dots,P_{n-1}$. Then the number of vertices in $\mathcal{I}$ counting multiplicities is
\begin{equation} \label{CurveVertexCount}
\sum_{A \in \mathcal{I}^{(0)}} m_A =  \textup{MV}_n (P_1,\dots,P_{n-1},P_1+\dots +P_{n-1}) \ .
\end{equation}
\end{cor}
\begin{remark}
Corollary~\ref{Cor:TransvIntVertexCount} generalizes \cite[Theorem 3.3]{Vigeland} where each $P_i$ is a standard simplex of the form $\textup{conv}\{ s_i\cdot e^{(i)} \cup\{0\} \,:\, 1\leq i\leq n \}$ where $e^{(i)}$ denotes the $i$-th unit vector and $s_i \in\Z_{>0}$. In this case \eqref{CurveVertexCount} gives $s_1\cdots s_{n-1}\cdot(s_1+\dots+s_{n-1})$ as the number of vertices counting multiplicities.
\end{remark}
\begin{proof}
For $k=n-1$ the sum in (\ref{VertexCount}) runs over all cells of type $(2,1,\dots,1)$, $(1,2,1,\dots,1)$, $\dots$, $(1,\dots,1,2)$. Hence,
\begin{align*}
\sum_{A \in \mathcal{I}^{(0)}} m_A &= \textup{MV}'_n (P_1,2;P_2,1;\dots,P_{n-1},1)+\dots+\textup{MV}'_n (P_1,1;P_2,1;\dots,P_{n-1},2).
\end{align*}
If the lattice $\sum_{i=1}^{n-1} \Z(P_i)$ spanned by the integer points of 
$P_1, \ldots, P_{n-1}$ 
coincides with $\Z^n$ then the volume forms $\textup{MV}_n$
and $\textup{MV}'_n$ coincide, and by the symmetry and linearity of the 
mixed volume \eqref{linearity} we get~\eqref{CurveVertexCount}.

In the general situation, by \cite{kkmd} we can multiply all polytopes
by a suitable factor $N$ to obtain a lattice $\sum_{i=1}^{n-1} \Z(N P_i)$
which coincides with $\Z^n$.
Then for all $l \in \N$, the lattice $\sum_{i=1}^{n-1} \Z(lNP_i)$ as well 
coincides with $\Z^n$, and hence~\eqref{CurveVertexCount} holds.
Since then both sides of~\eqref{CurveVertexCount} are polynomials
in $l$ and $N$, the equation thus follows with regard to the lattice
$\sum_{i=1}^{n-1} \Z(P_i)$ as well.
\end{proof}
We can also prove Corollary \ref{Cor:TransvIntVertexCount} independently of the dual approach by using stable intersections.
\begin{proof}
Define $\mathcal{J}:= X(f_1\odot \dots \odot f_{n-1})=X(f_1)\cup\dots\cup X(f_{n-1})$. We know that $\underbrace{\mathcal{J}\cap_{\textup{st}} \dots   \cap_{\textup{st}}\mathcal{J}}_{n\text{-times}}$ $= \mathcal{J}^{(0)}$. Since $\mathcal{I}\subset \mathcal{J}^{(1)}$ holds, this implies that $\mathcal{I}\cap_{\textup{st}} \mathcal{J}\subset \mathcal{J}^{(0)}$. Furthermore we have $\mathcal{I}\cap_{\textup{st}} \mathcal{J}\subset \mathcal{I}\cap\mathcal{J}=\mathcal{I}$ and $\mathcal{J}^{(0)}\cap\mathcal{I}=\mathcal{I}^{(0)}$ such that 
\[
\mathcal{I}^{(0)}= \mathcal{I}\cap_{\textup{st}}\mathcal{J} \ .
\]
The Newton polytope of $f_1\odot \dots \odot f_{n-1}$ is $P_1+\dots+P_{n-1}$.  Now using the tropical Bernstein theorem for stable intersections (Theorem \ref{theorem:TropicalBernstein}) we have that the number of points in $\mathcal{I}^{(0)}$ counted with multiplicities is $\textup{MV}_n (P_1,\dots,P_{n-1},P_1+\dots +P_{n-1})$.
\end{proof}

With similar techniques we will count now the number of unbounded faces in $\mathcal{I}=X_1\cap\dots\cap X_k$. Again, we formulate the result in a general manner though our main interest will later be the case $k=n-1$ and $j=1$, i.e. the number of unbounded edges in a tropical intersection curve.  
\begin{theorem} \label{theorem:unbounded}
The number of unbounded $j$-faces in $\mathcal{I}$ is 
\begin{equation}\label{UnboundedRays}
\sum_{\begin{smallmatrix}F=(P_1)^{v}+\dots+(P_{k})^{v} 
\end{smallmatrix}} 
\textup{MV}'_{n-j}((P_1)^{v},\dots,(P_{k})^{v})\ .
\end{equation}
Here the sum is taken over all $(n-j)$-faces $F$ of $P:=P_1+\dots+P_{k}$, $v\in\mathbb{S}^n$ is the outer unit normal vector of $F$ and $\textup{MV}'_{n-j}$ denotes the $(n-j)$-dimensional mixed volume taken with respect to the lattice defined by the face $F$.
\end{theorem}
\begin{proof}
As seen in Section \ref{section:duality} the unbounded $j$-faces of the union $X_1 \cup \dots \cup X_{k}$ correspond to $(n-j)$-dimensional cells in the boundary of $P=P_1+\dots+P_k$. So to count the unbounded $j$-faces in the intersection $\mathcal{I}$ we count mixed cells in all $(n-j)$-faces of $P$. Each face $F$ of $P$ has an outer unit normal vector $v$ and $F=(P_1)^v+\dots+(P_{k})^v$ where $(P_i)^v$ denotes the face of $P_i$ which is maximal with respect to $v$. 
So the number of unbounded $j$-faces counted with multiplicity (see Definition \ref{Def:Multiplicity})  which are dual to cells in $F$ is $\textup{MV}_{n-j}'((P_1)^{v},\dots,(P_{k})^{v})$ and the result follows. 
\end{proof}

\section{Lattice points and the genus of tropical intersection curves} \label{section:LatticePointsAndGenus}
Suppose an intersection curve of $n-1$ smooth tropical hypersurfaces in $\R^n$ is given. 
We apply the results of the last section to express its tropical genus $g$ (as defined below) in terms of mixed volumes of the Newton polyhedra corresponding to the defining hypersurfaces. Our goal is to prove that this genus coincides with the genus $\bar{g}$ of a toric variety $X$ that was obtained using the same Newton polytopes. Due to a result by Khovanskii \cite{Khovanskii78} the toric genus can be expressed via alternating sums of interior integer point numbers. To show that the combinatorial expressions for $g$ and $\bar{g}$ are equal we develop some aspects of a mixed Ehrhart theory. See \cite{katz-payne-2008} for 
related methods in the context of toric varieties.

\subsection{The genus via mixed volumes}
Assume in the following that the intersection curve $\mathcal{I}$ is connected and was obtained by a transversal intersection of $n-1$ hypersurfaces $X_1\cap\dots\cap X_{n-1}$ with Newton polytopes  $P_1,\dots,P_{n-1}$.
For such an intersection curve $\mathcal{I}$ in $\R^n$ define the \emph{genus} $g=g(\mathcal{I})$ as the number of independent cycles of $ \mathcal{I}$, i.e. its first Betti number.

Since $\mathcal{I}$ is a transversal intersection each vertex $A$ in $\mathcal{I}$ is dual to a cell $C$ of type $(1,\dots,1,2,1,\dots,1)$ in the privileged subdivision of $P_1+\dots+P_{n-1}$ . So $C$ is a sum of $n-1$ edges and one $2$-dimensional face $F_i$ of $P_i$. The \emph{degree} (or \emph{valence}) of $A$ is the number of outgoing edges (bounded and unbounded) in $\mathcal{I}$. Each such outgoing edge $A'$ is dual to an $(n-1)$-dimensional mixed cell $C'$ which is a facet of $C$. Hence the degree of $A$ equals the number of edges of the $2$-dimensional face $F_i$. 

Vigeland gave in \cite{Vigeland} an expression for the genus of a $3$-valent curve in terms of inner vertices and outgoing edges. The proof does not apply tropical properties of $\mathcal{I}$ and works for any $3$-valent graph with unbounded edges. Note that the vertices and edges are not counted with multiplicities in this statement. 
\begin{prop}[see \cite{Vigeland}] \label{lemma:genus}
For a $3$-valent tropical intersection curve $\mathcal{I}$ we have 
\[
2g-2= \#\{\text{vertices in }\mathcal{I}\} - \#\{\text{unbounded edges in }\mathcal{I}\}\ .
\]
\end{prop}
A tropical hypersurface $X_i$ is called \emph{smooth} if the maximal cells of its privileged subdivision are simplices of volume $\frac{1}{n!}$. If $\mathcal{I}$ is obtained as an intersection of smooth hypersurfaces then $\mathcal{I}$ is $3$-valent and each vertex and unbounded edge has multiplicity $1$.   

\begin{theorem}\label{Theorem:Genus}
Let $\mathcal{I}$ be a connected transversal intersection of $n-1$ smooth tropical hypersurfaces in $\R^{n}$ with Newton polytopes $P_1,\dots,P_{n-1}$. Then the genus $g$ of $\mathcal{I}$ is given by
\begin{equation}\label{eq:genus}
2g-2 = \textup{MV}_{n}\left(P_1,\dots,P_{n-1},\sum_{i=1}^{n-1} P_i\right) - \sum_{v}
\textup{MV}_{n-1}'((P_1)^{v},\dots,(P_{n-1})^{v})
\end{equation}
where $v$ runs over all outer unit normal vectors of $P_1+\dots+P_{n-1}$.
\end{theorem}
\begin{remark}
If the smoothness condition of the hypersurfaces $X_i$ is dropped the right hand side of \eqref{eq:genus} still gives an upper bound for $2g-2$.
\end{remark}
\begin{proof}
Using Corollary \ref{Cor:TransvIntVertexCount}, Theorem \ref{theorem:unbounded} and Proposition \ref{lemma:genus} we immediately get the result.
\end{proof}
In particular we see that under the conditions of Theorem~\ref{Theorem:Genus}
the genus only depends on the Newton polytopes $P_1,\dots,P_{n-1}$ and
we will write $g(P_1,\dots,P_{n-1})$ to denote this value.

\begin{example}\label{example:pick}
Consider this theorem in the case $n=2$. Here we just have one smooth tropical hypersurface $X_1$ with corresponding Newton polytope $P_1$. The genus $g$ of this curve equals the number of interior integer points of $P_1$. So Theorem \ref{Theorem:Genus} states that
\begin{align*}
2 \cdot \# \left\{ \begin{smallmatrix}\text{interior integer}\\ \text{points of }P_1\end{smallmatrix}\right\} -2
&= \textup{MV}_2(P_1,P_1)- \sum_{v \in\Sph^2}
\textup{MV}_{1}'((P_1)^{v})\\
&= 2\cdot \textup{vol}_2(P_1) - \# \left\{\begin{smallmatrix}\text{integer points on}\\ \text{the facets of }P_1\end{smallmatrix}\right\}\ .
\end{align*}
Hence Theorem \ref{Theorem:Genus} implies that 
\[
\text{vol}_2(P_1)= \# \left\{\begin{smallmatrix}\text{interior integer}\\ \text{points of }P_1\end{smallmatrix}\right\}
+ \frac{1}{2}\cdot \# \left\{\begin{smallmatrix}\text{integer points on}\\ \text{the facets of }P_1\end{smallmatrix}\right\} -1 \ 
\]
which is known as Pick's theorem for convex polygons (see \cite{TheBook}).
\end{example}

\subsection{The toric genus}
In \cite{Khovanskii78}, Khovanskii gave a formula for the genus of a complete intersection
in a toric variety. Let the variety $X$ in $(\C^*)^n$ be defined by a  non-degenerate system of equations $f_1=\dots=f_{k}=0$ with Newton polyhedra $P_1,\dots,P_{k}$ where each has full dimension $n$. Let $\bar{X}$ be the closure of $X$ in a sufficiently complete projective toric compactification. 
\begin{prop}[Khovanskii  \cite{Khovanskii78}] \label{the:Khovanskii}
If $\bar{X}$ is connected and has no holomorphic forms  of intermediate dimension, then the geometric genus $\bar{g}$ of $X$  can be calculated by the formula
\begin{equation}\label{eq:Khovanskii}
\sum_{\emptyset\neq J\subset [k]}(-1)^{k-|J|} \cdot B^{+}(\sum_{j\in J} P_j)
\end{equation}
where $B^{+}(P)$ denotes the number of interior integer points of the lattice polytope $P$ and $[k]:=\{1,\dots,k\}$.
\end{prop}
Thus for any variety satisfying the conditions of Proposition \ref{the:Khovanskii}, the genus only depends on $P_1,\dots,P_k$. We call this value $\bar{g}(P_1,\dots,P_k)$.

We are ready now to state our theorem comparing the genus of tropical and toric intersection curves. 
\begin{theorem}\label{the:CompareGenus}
Let $P_1,\dots,P_{n-1}\subset\R^n$ be full-dimensional lattice polytopes. Then the tropical and the toric genus with respect to $P_1,\dots,P_{n-1}$ coincide, i.e.
\[
\bar{g}(P_1,\dots,P_{n-1}) = g(P_1,\dots,P_{n-1})\ .
\]
\end{theorem}
We will prove this theorem by showing that the combinatorial quantities of Proposition \ref{the:Khovanskii} and Theorem \ref{Theorem:Genus} are the same, i.e.
\begin{align}
 &\frac{1}{2}\,\textup{MV}_{n}(P_1,\dots,P_{n-1},\sum_{i=1}^{n-1} P_i) - \frac{1}{2}\sum_{v\in\Sph^n}
 \textup{MV}_{n-1}'((P_1)^{v},\dots,(P_{n-1})^{v}) +1 \nonumber \\
 =\ &\sum_{\emptyset\neq J\subset [n-1]}(-1)^{n-1-|J|}\cdot  B^{+}(\sum_J P_j)\ . \label{eq:TheEquation}
\end{align}

That \eqref{eq:TheEquation} holds for $n=2$ can be seen in Example \ref{example:pick} when Pick's theorem is assumed to be given. For higher dimensions further combinatorial considerations are necessary.

\subsection{The unmixed case of Theorem \ref{the:CompareGenus}} 
We study the surface volume and the number of integer points of a \emph{lattice complex}, i.e. a bounded polyhedral complex with vertices in $\Z^n$.
By using a classical result of Macdonald we will be able to employ this to prove the \emph{unmixed case}  $P_1=\dots=P_{n-1}$ of \eqref{eq:TheEquation}. 

Let $\chi(Q)$ denote the Euler-Poincar\'e characteristic of a polyhedral complex $Q$. For simplicity we set $B(0\cdot Q):=\chi(Q)$ and by $\partial Q$ we denote the boundary complex of $Q$.

\begin{theorem}\label{the:PickOnSurface}
Let $Q$ be a pure $n$-dimensional lattice complex. Then 
\begin{equation}\label{eq:n-dimUnmixed}
\sum_{F \textup{ facet of }Q}(n-1)!\, \textup{vol}_{n-1}'(F) =  (-1)^{n-1}\sum_{k=0}^{n-1}(-1)^{k}\binom{n-1}{k} B(k\cdot \partial Q) \ .
\end{equation}
\end{theorem}
\begin{proof}
We first consider the case where the facets of $k\cdot Q$ admit a
unimodular triangulation (i.e., a triangulation into
simplices $\Delta$ of volume $\frac{1}{(\dim \Delta) ! }$) 
with respect to the lattices defined by the facets. 
Let $f_i$ be the number of $i$-dimensional faces of this simplicial complex. Note that the left hand side of \eqref{eq:n-dimUnmixed} counts the number of $(n-1)$-dimensional faces, i.e.
\[
f_{n-1}= \sum_{F \text{ facet of }Q}(n-1)!\, \textup{vol}_{n-1}'(F)\ .
\]
Each $i$-dimensional face of our complex is a fundamental lattice simplex. The number of interior integer points of a fundamental lattice simplex $\Delta$ of dimension $i$ stretched by a factor of $k\geq 1$ is equal to 
\[
\# \left\{x\in\N^i\,:\, x_j\geq 1\text{ and } \sum_j x_j \leq k-1   \right\} = \binom{k-1}{i}\ . 
\]
Hence we have for $k\geq 1$ that
\[
B(k\cdot \partial Q) = \sum_{i=0}^{k-1} \binom{k-1}{i} f_i \ .
\]
Up to the term for $k=0$ 
the sum on the right hand side of \eqref{eq:n-dimUnmixed} evaluates to
\begin{align}
&\sum_{k=1}^{n-1}(-1)^{n-1-k}\binom{n-1}{k}\sum_{i=0}^{k-1} \binom{k-1}{i} f_i \nonumber \\
=\ & (-1)^{n-1}\sum_{i=0}^{n-2} f_i \sum_{k=1+i}^{n-1} (-1)^k \binom{n-1}{k}\binom{k-1}{i}\nonumber \\
=\ & \sum_{i=0}^{n-2} f_i \sum_{r=0}^{n-2-i} (-1)^r \binom{n-1}{r}\binom{n-2-r}{n-2-i-r}
 \label{eq:surf2} 
\end{align}
where we substituted $r=n-k-1$ to obtain the last equation.
Using the following binomial identity (see e.g. \cite[p. 149]{Gruenbaum})
\[
\text{For }0\leq c\leq a: \quad \sum_{i=0}^c (-1)^i \binom{b}{i}\binom{a-i}{c-i} = \binom{a-b}{c}
\]
yields that the right hand side in \eqref{eq:n-dimUnmixed} equals
\[
(-1)^{n-1} \chi (\partial Q) + \sum_{i=0}^{n-2} f_i \binom{-1}{n-2-i} = (-1)^{n-1} \chi (\partial Q)+ \sum_{i=0}^{n-2} (-1)^{n-2-i} f_i\ .
\]
By the Euler-Poincar\'e formula  
$
\chi(\partial Q)= \sum_{i=0}^{n-1} (-1)^i f_i  
$ (see \cite{Bredon})
this expression simplifies to $f_{n-1}$ which
proves the theorem for the unimodular case.

In the general situation, we can proceed as in the proof of
Corollary~\ref{Cor:TransvIntVertexCount}. By \cite{kkmd},
there exists an $N \in \N$ such that $NQ$ admits a unimodular triangulation.
Then for all $l \in \N$, $l N Q$ as well admits a unimodular
triangulation and~\eqref{eq:n-dimUnmixed} holds.
Since both sides of the equation are polynomials
in $l$ and $N$, the equation follows for $Q$ as well.
\end{proof}

By combining Theorem \ref{the:PickOnSurface} and the generalization of Pick's theorem by Macdonald (see Reeve \cite{Reeve57} for the $3$-dimensional case) we get the unmixed version of \eqref{eq:TheEquation}.
\begin{prop}[Macdonald \cite{MacDonald62}]\label{the:MacDonald}
 Let $P$ be a pure $n$-dimensional lattice complex. 
Then with the notations from above we have
\[
\frac{n-1}{2} n!\, \textup{vol}_n(P) = \sum_{k=0}^{n-1} (-1)^{n-1-k}\binom{n-1}{k} \left[B(k\cdot P) - \frac{1}{2} B(k\cdot \partial P)\right] 
\ .
\]
\end{prop}
\begin{cor}\label{the:n-dimUnmixed}
For $n$-dimensional lattice polytopes $P$ we have
\[
\frac{n-1}{2} n!\, \textup{vol}_n(P)-\frac{1}{2} \sum_{F \text{ facet of }P}(n-1)!\, \textup{vol}_{n-1}'(F)  \ +1 
= \sum_{k=1}^{n-1}(-1)^{n-1-k}\binom{n-1}{k}B^{+}(k\cdot P)\ .
\]
\end{cor}
\begin{proof}
Knowing Proposition \ref{the:MacDonald} and the fact that $B(P)-\frac{1}{2}B(\partial P)=B^+(P)+\frac{1}{2}B(\partial P)$ we still have to show that
\begin{align*}
&(-1)^{n-1} \left[ \sum_{k=1}^{n-1} (-1)^{k}\binom{n-1}{k} B(k\cdot \partial P)\right] + 2(-1)^{n-1}(\chi(P)-\chi(\partial P))\\
=\ &\sum_{F \text{ facet of }P} (n-1)!\, \textup{vol}_{n-1}'(F)\ -2 \, .
\end{align*}
Since the Euler-Poincar\'e formula implies $\chi(P)= \chi(\partial P) +(-1)^n$,  the last equation reduces to the statement of Theorem \ref{the:PickOnSurface}.
\end{proof}

\subsection{The mixed case of Theorem \ref{the:CompareGenus}}
 For lattice polytopes $P_1,\dots,P_{k}\subset\R^{n}$ and $t\in\N$ consider the following version of a \emph{mixed Ehrhart polynomial},
\[
\ME_{P_1,\dots,P_k}(t) := \sum_{ \emptyset \neq J \subset \{1,\dots,k\}} (-1)^{k-|J|} B(t\cdot \sum_{j\in J} P_j)\, .
\]
As we will see this alternating sum of Ehrhart polynomials turns out to have a very simple structure. Namely all coefficients of $t^r$ for $1\leq r<k$ vanish and in the case of our main interest ($k=n-1$) the remaining coefficients have a nice interpretation in terms of mixed volumes. These results will prove
 \eqref{eq:TheEquation} in the mixed case.

Clearly $\ME_{P_1,\ldots,P_{k}}(t)$ is a polynomial in $t$ of degree at most $n$ since it is the alternating sum of Ehrhart polynomials (see Section \ref{section:Ehrhart}):
\begin{align*}
\ME_{P_1,\dots,P_k}(t) &= \sum_{ \emptyset \neq J \subset [k]} (-1)^{k-|J|} E_{\sum_J P_j}(t) \\
  &= \sum_{r=0}^n t^r \left( \sum_{ \emptyset \neq J \subset [k]} (-1)^{k-|J|} e_r(\sum_J P_j) \right)\ .
\end{align*}
We denote the coefficients of this polynomial by $\me_r(P_1,\dots,P_{k})$. To prove Theorem \ref{the:CompareGenus} we have to consider the alternating sum of expressions in $B^+$ rather then in $B$, i.e. we deal with interior integer points instead of just integer points. Fortunately the Ehrhart reciprocity \eqref{eq:EhrhartReciprocity} allows to translate each result in terms of $B$ to results in terms of $B^+$. Namely we have that
\begin{equation}\label{eq:MEforB+}
\sum_{ \emptyset \neq J \subset [k]} (-1)^{k-|J|} B^+(t\cdot \sum_{j\in J} P_j)= \sum_{r=0}^n t^r \cdot (-1)^{n+r} \cdot \me_r(P_1,\dots,P_{k})\ .
\end{equation}

As indicated above the following key lemma holds.
\begin{lemma}\label{lemma:CoefficientsVanish} 
For any polytopes $P_1,\dots,P_k$ and any $\Z^n$-valuation $\varphi_r$ which is homogeneous of degree $r \in \{1, \ldots, k-1\}$ we have that
\[
\sum_{\emptyset \neq J \subset [k]}(-1)^{k-|J|} \varphi_r(\sum_J P_j) = 0\ .
\]
\end{lemma}
In particular this implies that $\me_r(P_1,\ldots,P_{k}) = 0$ for $1\leq r <k$.
\begin{proof}
By McMullen's result on homogeneous valuations (see Proposition \ref{lemma:McMullen}) we obtain
\begin{equation*}
\sum_{\emptyset \neq J \subset [k]}(-1)^{k-|J|} \varphi_r(\sum_J P_j) =\sum_{\emptyset \neq J \subset [k]}(-1)^{k-|J|} \sum_{r_1,\dots,r_{|J|}}\binom{r}{r_1 \dots r_{|J|}} \varphi'_r( P_{j_1},r_1;\dots;P_{j_{|J|}},r_{|J|}) \ .
\end{equation*}
Here the $\varphi'_r( P_{j_1},r_1;\dots;P_{j_{|J|}},r_{|J|})$ are mixed valuations which we do not need to state more explicit. We write the right hand side of the previous equation slightly different as
\begin{equation}\label{ProofEq:3}
(-1)^{n-1} \hspace{-0.3cm} \sum_{\emptyset \neq J \subset [k]} \hspace{-0.3cm} (-1)^{|J|} \hspace{-0.5cm} \sum_{\begin{smallmatrix} s_1,\dots,s_{k}\geq 0\\\sum s_i = r \\ s_i=0 \text{ if } i\in[k]\backslash J \end{smallmatrix}}  \hspace{-0.3cm} \binom{r}{s_1 \dots s_{n-1}} \varphi'_r( P_1,s_1;\dots;P_{k},s_{k}) \ .
\end{equation}
Now fix $s_1,\dots,s_{k}\geq 0$ and ask for which sets $J$ does $\varphi'_r( P_1,s_1;\dots;P_{k},s_{k})$ appear in the inner sum of \eqref{ProofEq:3}. Denote by $J_s$ the set of indices $i$ for which $s_i \neq 0$ then $\varphi'_r( P_1,s_1;\dots;P_{k},s_{k})$ appears whenever $J\supset J_s$. If it appears then always with the same multinomial coefficient but possibly different sign depending on the number of elements in $J$. 

Let $\alpha_s$ be the number of elements in $[k]$ which are not in $J_s$. Then we can write 
\begin{equation}\label{ProofEq:4}
 \sum_{\emptyset \neq J \subset [k]}(-1)^{k-|J|} \varphi_r(\sum_J P_j) =(-1)^{k} \sum_{\begin{smallmatrix} s_1,..,s_{k}\geq 0\\\sum s_i = r \end{smallmatrix}} 
A(s)\cdot \varphi'_r( P_1,s_1;\dots;P_{k},s_{k})
\end{equation}
where
\[
A(s)= (-1)^{|J_s|} \binom{r}{s_1 \dots s_{k}} \sum_{i=0}^{\alpha_s} (-1)^i \binom{\alpha_s}{i} \ .
\]
Now $\sum_{i=0}^{\alpha_s} (-1)^i \binom{\alpha_s}{i}$ equals $0$ if $\alpha_s > 0$ and $1$ if $\alpha_s=0$. Since $1 \le r<k$ the case $\alpha_s=0$ can not occur and hence \eqref{ProofEq:4} vanishes for $1\leq r<k$. 
\end{proof}

We will now restrict ourselves to the case $k=n-1$. According to Lemma~\ref{lemma:CoefficientsVanish}  most coefficients of $\ME_{P_1,\dots,P_{n-1}}(t)$ vanish and the following theorem will show that the remaining coefficients can be explicitly expressed using classical results on Ehrhart coefficients (see \eqref{eq:CoeffIdentities}).
\begin{theorem}\label{the:MixedEhrhart}
$\ME_{P_1,\dots,P_{n-1}}(t)$ is a polynomial in $t$ of degree $n$ and we have
\begin{eqnarray}
\ME_{P_1,\dots,P_{n-1}}(t) &=&t^n\cdot \frac{1}{2}\textup{MV}_{n}\left(P_1,\dots,P_{n-1},\sum_{i=1}^{n-1} P_i\right) \nonumber\\
&\quad& + t^{n-1}\cdot \frac{1}{2}\sum_{v\in\Sph^n}
\textup{MV}_{n-1}'((P_1)^{v},\dots,(P_{n-1})^{v}) \ +\ (-1)^{n}\ . \nonumber 
\end{eqnarray}
\end{theorem}
From this follows equation \eqref{eq:TheEquation} and therefore Theorem \ref{the:CompareGenus} by setting $t=1$ and using the Ehrhart reciprocity as in \eqref{eq:MEforB+}.
\begin{proof}
By Lemma \ref{lemma:CoefficientsVanish}, it remains to compute the coefficients of $t^n$, $t^{n-1}$ and $t^0$. 
The constant coefficient $e_0(P)$ equals $1$ for every $P$ (see \eqref{eq:CoeffIdentities}). Hence
\begin{align*}
\me_0(P_1,\dots,P_{n-1})&= \sum_{\emptyset \neq J \subset [n-1]}(-1)^{n-1-|J|}\cdot 1\\
& = (-1)^{n-1} \sum_{\emptyset \neq J \subset [n-1]}(-1)^{|J|}= (-1)^n \ .
\end{align*}
Considering $e_n(P)$, \eqref{eq:CoeffIdentities} states that  the coefficient of $t^n$ in the Ehrhart polynomial equals $\vol_n(P)$. Hence 
\begin{equation}\label{ProofEq:5}
\me_n(P_1,\dots,P_{n-1})=(-1)^{n-1}\sum_{1\leq i_1<\dots<i_k \leq n-1}(-1)^{k} \vol_n(P_{i_1}+\dots+P_{i_k})\ .
\end{equation}
Basic results on mixed volumes (see, e.g. \cite[\S 7.4]{CLO2}) show 
\[
\vol_n(P_{i_1}+\dots+P_{i_k})= \sum_{\begin{smallmatrix}j_1+\dots+j_k=n\\ j_s\geq 0 \end{smallmatrix}} \frac{1}{j_1!\cdots j_k!} \mv_n(P_{i_1},j_1;\dots;P_{i_k},j_k)\ ;
\]
thus the right hand side of \eqref{ProofEq:5} can be written as
\begin{align}
&(-1)^{n-1}\sum_{1\leq i_1<\dots<i_k \leq n-1}(-1)^{k} \sum_{\begin{smallmatrix}j_1+\dots+j_k=n\\ j_s\geq 0 \end{smallmatrix}} \frac{1}{j_1!\cdots j_k!} \mv_n(P_{i_1},j_1;\dots;P_{i_k},j_k) \nonumber \\ 
=\ &(-1)^{n-1}\sum_{\emptyset \neq J \subset [n-1]}(-1)^{\# J} \sum_{\begin{smallmatrix} s_i\geq 0,\  \sum s_i = n \\ s_i=0 \text{ if } i\in[n-1]\backslash J \end{smallmatrix}}\frac{\mv_n(P_1,s_1;\dots;P_{n-1},s_{n-1})}{s_1!\cdots s_{n-1}!} \label{ProofEq:7} 
\end{align}
With the same notation $J_s$ and $\alpha_s$ as in the proof of Lemma \ref{lemma:CoefficientsVanish} we see that $\mv_n(P_1,s_1;$ $\ldots ; P_{n-1},s_{n-1})$ appears in the inner sum of \eqref{ProofEq:7} whenever $J_s\subset J$.  Using this in \eqref{ProofEq:7} we get
\begin{equation}\label{ProofEq:8}
\me_n(P_1,\dots,P_{n-1})= (-1)^{n-1} \sum_{s_i\geq 0,\  \sum s_i = n} A'(s)\cdot \mv_n(P_1,s_1;\dots;P_{n-1},s_{n-1})
\end{equation}
where $A'(s)=\frac{(-1)^{|J_s|}}{ s_1!\cdots s_{n-1}!}\sum_{i=0}^{\alpha_s} (-1)^i \binom{\alpha_s}{i}$.
As seen before $A'(s)=0$ for $\alpha_s\neq 0$. Hence only terms with  $\alpha_s=0$ (i.e. $J_s=[n-1]$) remain in which case $\sum_{i=0}^{\alpha_s} (-1)^i \binom{\alpha_s}{i}=1$ and we obtain
\begin{align*}
\me_n(P_1,\dots,P_{n-1}) &=(-1)^{n-1} \sum_{s_i\geq 1,\  \sum s_i = n}(-1)^{n-1} \frac{\mv_n(P_1,s_1;\dots;P_{n-1},s_{n-1})}{2!\cdot 1! \cdots 1!} \\
&=\frac{1}{2}\textup{MV}_{n}\left(P_1,\dots,P_{n-1},\sum_{i=1}^{n-1} P_i\right)\ .
\end{align*}

The coefficient of $t^{n-1}$ can be computed using the same combinatorial trick as before only that we start here with $e_{n-1}(P)=\frac{1}{2}\sum_{F \text{ facet of }P}\vol_{n-1}'(F) $ (see \eqref{eq:CoeffIdentities}). 
\begin{align}
\me_{n-1}(P_1,\dots,P_{n-1}) 
&= (-1)^{n-1}\sum_{1\leq i_1<\dots<i_k \leq n-1}(-1)^{k}\,\frac{1}{2} \sum_{\begin{smallmatrix}F \text{ facet of }\\ P_{i_1}+\dots+P_{i_k}\end{smallmatrix}}\vol_{n-1}'(F) \nonumber \\
&=(-1)^{n-1}\sum_{1\leq i_1<\dots<i_k \leq n-1}(-1)^{k}\, \frac{1}{2}\sum_{ v \in \Sph^n}  \vol_{n-1}'\left((P_{i_1}+\dots+P_{i_k})^v \right) \ ,    \label{ProofEq:10}
\end{align}
where the last equation holds since $\vol_{n-1}'\left((P_{i_1}+\dots+P_{i_k})^v \right) $ vanishes whenever $v$ is not a facet normal of $P_{i_1}+\dots+P_{i_k}$. 
Since $(P_{i_1}+\dots+P_{i_k})^v= (P_{i_1})^v+\dots+(P_{i_k})^v$ holds \eqref{ProofEq:10} can be written as
\begin{equation}\label{ProofEq:11}
\frac{1}{2}\sum_{v \in \Sph^n} \left[ (-1)^{n-1}\sum_{1\leq i_1<\dots<i_k \leq n-1}(-1)^{k}\vol_{n-1}'\left((P_{i_1})^v+\dots+(P_{i_k})^v \right) \right]\ .
\end{equation}
With the same method as before (starting from equation \eqref{ProofEq:5}) we can show that the term in the large brackets in \eqref{ProofEq:11}
equals
$
\textup{MV}_{n-1}'((P_1)^{v},\dots,(P_{n-1})^{v})\ .
$
Now finally using this in \eqref{ProofEq:11} yields 
\[
\me_{n-1}(P_1,\dots,P_{n-1}) = \frac{1}{2}\sum_{v \in \Sph^n} \textup{MV}_{n-1}'((P_1)^{v},\dots,(P_{n-1})^{v})\ .
\]
This proves our theorem and therefore shows that the tropical and toric genus coincide.
\end{proof}
\smallskip

\noindent{\bf Acknowledgement:} 
We wish to thank Hartwig Bosse, Benjamin Nill, and Raman Sanyal
for helpful discussions on several aspects of the work described in this paper. The pictures in example \ref{Example1} were made using the tropical maple package of N.~Grigg, \texttt{http://math.byu.edu/tropical/maple}.


\begin{thebibliography}{}

\bibitem{TheBook} M.~Aigner and G.~Ziegler. {\it Proofs from the Book.} Third Edition. (Springer-Verlag, Berlin), 2004. 

\bibitem{Barvinokbook} A.~Barvinok. {\it Integer Points in Polyhedra.} (Z\"urich lectures in Adv. Math., Europ. Math. Soc.), 2008. 

\bibitem{Bernstein} D.N. Bernstein. The number of roots of a system of equations. {\it Funkcional. Anal. i Prilo\v zen.} 9(3):1--4, 1975.

\bibitem{BertrandBihan} B.~Bertrand and  F.~Bihan. Euler characteristic of real nondegenerate tropical complete intersections. Preprint, \texttt{arXiv:math/0710.1222}, 2007. 

\bibitem{Betke85} U.~Betke and M.~Kneser. Zerlegungen und Bewertungen von Gitterpolytopen. {\it J. Reine Angew. Math.} 358:202--208, 1985.

\bibitem{Bredon} G.E.~Bredon. {\it Topology and Geometry.} (Graduate Texts in Mathematics, Vol. 139. Springer, New York), 1993.

\bibitem{CLO2} D.A.~Cox, J.~Little and D.~O'Shea. {\it Using Algebraic Geometry.} Second Edition. (Graduate Texts in Mathematics, Vol. 185. Springer, New York), 2005.

\bibitem{SturmfelsFeichtner} A.~Dickenstein,  E.M.~Feichtner and B.~Sturmfels. Tropical discriminants. {\it J. Amer. Math. Soc.} 20:1111--1133, 2007.

\bibitem{Draisma08} J.~Draisma. A tropical approach to secant dimensions. {\it J. Pure Appl. Algebra} 212(2):349--363, 2008.

\bibitem{Ehrhart67} E.~Ehrhart. Sur un probl\`eme de g\'eometrie diophantinenne lin\'eaire, I. {\it J. Reine Angew. Math.} 226:1--29, 1967.

\bibitem{EinsiedlerKapranov} M.~Einsiedler, M.~Kapranov and D.~Lind. Non-archimedean amoebas and tropical varieties. {\it J. Reine Angew. Math.} 601:139--157, 2006.

\bibitem{Huber} B.~Huber and B.~Sturmfels. A polyhedral method for solving sparse polynomial systems. {\it Math. Comp.} 64:1541--1555, 1995.

\bibitem{GathmannMarkwig} A.~Gathmann and H.~Markwig. The numbers of tropical plane curves through points in general position. {\it J. Reine Angew. Math.} 602:155--177, 2007.

\bibitem{Gruenbaum} B.~Gr\"unbaum. {\it Convex Polytopes.} Second Edition. (Graduate Texts in Mathematics, Vol. 221. Springer, New York), 2003.

\bibitem{Itenberg} I.~Itenberg, V.~Kharlamov and E.~Shustin.  Welschinger invariant and enumeration of real rational curves. {\it Int. Math. Res. Not.} 49:2639--2653, 2003.

\bibitem{Katz08} E.~Katz. A Tropical Toolkit. 
 \emph{Expo. Math.} 27:1--36, 2009.

\bibitem{kkmd}
G.~Kempf, F.F.\ Knudsen, D.~Mumford and B.~Saint-Donat,
\emph{Toroidal embeddings I.}
Lecture Notes in Mathematics, vol. 339,
Springer, Berlin, 1973.

\bibitem{Khovanskii78} A.~Khovanskii. Newton polyhedra and the genus of complete intersections. {\it Funkcional. Anal. i Prilo\v zen.} 12(1):51--61, 1978.

\bibitem{MacDonald62} I.~G.~Macdonald. The volume of a lattice polyhedron. {\it Proc. Camb. Phil. Soc.} 59:719--726, 1962.

\bibitem{McMullen77} P.~McMullen. Valuations and Euler-type relations on certain classes of convex polytopes. {\it Proc. London. Math. Soc.} (3)35:113--135, 1977. 

\bibitem{McMullen93} P.~McMullen. Valuations and dissections. In: P.M.~Gruber and J.M.~Wills, editors, {\it Handbook of Convex Geometry,} volume B. North Holland, Amsterdam, 933--988, 1993.

\bibitem{Mikhalkin04b} G.~Mikhalkin. Enumerative Tropical Algebraic Geometry in $\R^2$. {\it J. Amer. Math. Soc.} 18:313--377, 2005.

\bibitem{MikhalkinPairsofPants} G.~Mikhalkin. Decomposition into pairs-of-pants for complex algebraic      hypersurfaces. {\it Topology}, 43:1035--1065, 2004.

\bibitem{Mikhalkin04} G.~Mikhalkin. Tropical geometry and its applications.  {\it International Congress of Mathematicians}, Vol. 2 (Eur. Math. Soc., Z\"urich), 827--852, 2006. 

\bibitem{SiebertNishinou} T.~Nishinou and B.~Siebert. Toric degenerations of toric varieties and tropical curves. {\it Duke Math. J.} 135:1--51, 2006.

\bibitem{katz-payne-2008}
E.~Katz and S.~Payne. 
Piecewise polynomials, Minkowski weights, and localization
  on toric varieties,
\emph{Algebra Number Theory} 2:135--155, 2008.

\bibitem{Reeve57} J.~E.~Reeve. On the volume of lattice polyhedra. {\it Proc. London. Math. Soc.} 34:378--395, 1959.

\bibitem{FirstStepsTropical} J.~Richter-Gebert, B.~Sturmfels and T.~Theobald. First steps in tropical geometry. In: {\it Proc. Idempotent mathematics and mathematical physics,}  Vienna 2003.  {\it Contemp. Math.} Vol. 377 (Amer. Math. Soc., Providence, RI), 289--317, 2005.

\bibitem{Schneider} R.~Schneider. {\it Convex bodies: The Brunn-Minkowski theory.} Encyclopedia of Mathematics and its Applications. (Cambridge University Press, Cambridge 44), 1993.

\bibitem{SpeyerTropicalLinearSpaces} D.~E.~Speyer. Tropical linear spaces.  {\it SIAM J. Discrete Math.}  22(4):1527--1558, 2008.

\bibitem{SpeyerSturmfels} D.~E.~Speyer and B.~Sturmfels. The tropical Grassmannian. {\it Adv. Geom.} 4:389--411, 2004.

\bibitem{sturmfels-tevelev-2008}
B.~Sturmfels and J.~Tevelev.
\newblock Elimination theory for tropical varieties.
\newblock {\em Math. Res. Lett.}, 15:543--562, 2008.

\bibitem{Vigeland} M.~D.~Vigeland. Tropical complete intersection curves. Preprint, \texttt{arXiv:math/0711.1962}, 2007.

\end{thebibliography}
\end{document}